\input AHTOH-E.STY
\UDC{
512.543.72 
+
512.544.33 
}
\MSC{
20F70,   
20F34,   
20F16    
}

\title{The Klein bottle group is not strongly verbally closed,
though awfully close to being so}

\author{%
Anton A. Klyachko
}
\address{
\myAddressW
}

\grants{\RFBR19-01-00591}

\abstract{%
According to Mazhuga's theorem,
the fundamental group $H$ of any
connected surface,
possibly except for the Klein bottle,
is
a retract of each
finitely generated group
containing $H$ as a verbally closed subgroup.
We prove that the Klein bottle group is indeed an exception
but has a very close property.
}

\s 0.
Introduction

A subgroup $H$ of a group $G$ is called \emph{verbally closed}
[MR14]
(%
see also
[Rom12],
[RKh13],
[Mazh17],
\hbox{[KlMa18]},
[KMM18],
[Mazh18],
[Bog18],
[Bog19],
[Mazh19],
[RT19]%
)
if
any equation of the form
$$
w(x,y,\dots)=h,
\qbox{where $w(x,y,\dots)$ is an element of the free group
$F(x,y,\dots)$ and $h\in H$,
}
$$
having a solution in $G$ has a solution in $H$.
If each finite system of equations with coefficient from~$H$
$$
\{w_1(x,y,\dots)=1, \dots, w_m(x,y,\dots)=1\},
\qbox{where $w_i\in H*F(x,y,\dots)$ (and $*$ stands for the
free product),}
$$
having a solution in $G$ has a solution in $H$, then $H$
is called \emph{algebraically closed} in $G$.

Algebraic closedness is a stronger property than verbal
closedness;
but these properties turn out to be
equivalent in many cases.
A group $H$ is called \emph{strongly verbally closed}
[Mazh18]
if it is algebraically closed
in any group containing $H$ as a verbally closed
subgroup. Thus, the verbal closedness is a
subgroup property, while the strong verbal
closedness is an abstract group property.
The class of strongly verbally closed groups is fairly wide,
see the papers mentioned
above. For instance, in [Mazh18], it was proven that
\disp{
\sl\narrower\narrower\narrower\narrower
the fundamental group of any
connected surface,
possibly except for the Klein bottle,
is
strongly verbally closed.
}%
The intriguing unique possible exception arose
as follows:
\-
almost all surface groups
are similar to free
groups in a sense; for such groups,
an ``industrial" method works;
this method goes back
to the very first paper [MR14] on this subject and
is
based on the use of the Lie words [Lee02]
or their analogues;
\-
the remaining several groups are either
\itemitem{}
abelian and, hence,
strongly verbally closed
(by a very simple reason) [Mazh18],
\itemitem{or}
this exceptional Klein bottle group
$K=\pres<a,b|a^b=a^{-1}>$,
which is neither free-like nor abelian.

\enditem
We obtain a
natural addition to Mazhuga's theorem:
\disp{\sl
the fundamental group of the Klein bottle is
not strongly verbally closed.
}%

A similar situation appeared
some time ago:
in [KlMa18], it was proven a general theorem
with a unique possible exception:
\disp{\sl
all non-dihedral virtually free groups
containing no non-identity finite normal subgroups
are strongly verbally closed
}%
(and the condition of absence of finite normal subgroups
cannot be removed);
later, it turned out [KMM18] that
\disp{\sl
the infinite dihedral group is
strongly verbally closed too,
}%
and this humble result
was more difficult
than the above mentioned general theorem
obtained in \hbox{[KlMa18]} by the ``industrial"\ method.
Generally, if the reader
takes some particular
his favorite nonabelian group,
far from free (e.g., a finite group),
then it would likely be difficult to decide
whether
this group is strongly verbally
closed:
{neither}
proving strong verbal closedness is easy,
{nor}
disproving this property is easy
(actually,
it is not too easy even to give an example of a
group
not being strongly verbally closed [KlMa18]).

The Klein bottle group
$K=\pres<a,b|a^b=a^{-1}>$
and the
infinite dihedral group
$D_\infty=\pres<a,b|a^b=a^{-1},\ b^2=1>$
are similar of course.
We use this similarity,
apply the result of [KMM18],
and conclude that the Klein bottle group,
though not strongly verbally closed, has a very similar property.

If a group~$H$ is \emph{equationally Noetherian} (i.e. any system
of
equations over $H$ with finitely many unknowns is equivalent
to its
finite subsystem),
then
the algebraic closedness is equivalent
to
the ``local retractness" [KMM18]:
\disp{
\sl\narrower\narrower\narrower
an equationally Noetherian
subgroup $H$
of a
group $G$
is algebraically closed in $G$ if and only if $H$ is
a retract
{\rm\(i.e. the image of an
endomorphism $\rho$ such that $\rho\o\rho=\rho$\)}
of each
finitely generated over $H$
subgroup of $G$ \rm(i.e. a subgroup of the form $\gp{H\cup X}$,
where $X\subseteq G$ is a finite set).
}%
All surface groups
are linear [New85];
and all linear groups are equationally Noetherian [BMR99].
Thus, our main
(and sole)
result can be stated as follows.

\Th.
The fundamental group $K$ of
the Klein bottle
{\rm(unlike all other surface groups)}
embeds into a
finitely generated group $G$
as a verbally closed subgroup that is not a retract of $G$.
However, any such $G$ has an index-two subgroup
containing $K$ as its retract.

In the next section, we give
an example proving the first assertion of the theorem
(i.e. that $K$ is not strongly verbally closed).
Section 2 contains
auxiliary lemmata.
In the last section, we prove the second assertion of the theorem.

The author expresses his deep gratitude to Veronika Miroshnichenko;
we have spent a long time
trying together to solve this problem
that seemed to be hard.
I also thank Andrey Mazhuga for reading a
draft of this text.

\smallskip
\noindent
{\bf Our notation}
is fairly standard. Note only that, if
$x$ and $y$ are elements of a group, then $x^y$
denotes $y^{-1}xy$,
The commutator~$[x,y]$ is $x^{-1}y^{-1}xy$.
If
$X$~is a subset of a group,
then
$\gp X$,
$\nc X$,
and $C_H(X)$
stand for
the subgroup generated by $X$,
the normal closure of $X$,
and
the centraliser of~$X$ in $H$
(where $H$~is a subgroup).
The symbol $\gp x_k$ denotes
the cyclic group of order $k$ generated by $x$.
The free group with a basis $x_1,\dots,x_n$ is denoted
as $F(x_1,\dots,x_n)$.

\s 1.
An example

Let
$V_4=\{1,d_1,d_2,d_3\}$ be the Klein four-group
(i.e. the noncyclic group of order four).
Consider the semidirect product
$$
G=\Bigl(V_4\times\gp{b}_\infty\Bigr)
\semitimes
\Bigl(
\gp{a_1}_\infty\times\gp{a_2}_\infty\times\gp{a_3}_\infty
\Bigr),
\qbox{where the action is }
a_i^b=a_i^{-1},\
a_i^{d_i}=a_i,\
a_i^{d_j}=a_i^{-1} \hbox{ for all }i\ne j.
$$

\Proposition.
The subgroup $K=\gp{a,b}\subset G$, where $a=a_1a_2a_3$,
\(isomorphic to the Klein bottle group\)
is
verbally closed in $G$, but not a retract.

\Proof
The subgroup $K$ is not a retract, because,
a hypothetical retraction $G\to K$ should map
finite-order
elements~$d_i$
to 1 (as $K$ is torsion-free);
then, the relation~$a_i^{d_j}=a_i^{-1}$ would show that the images
of $a_i$ are also of finite order and, hence, are 1 too;
therefore, the image of $a=a_1a_2a_3\in K$ is also 1 that
contradicts the fixedness of elements of $K$ under the retraction.

It remains to show that $K$ is verbally closed in $G$,
i.e.
any equation of the form
$$
w(x,y,\dots)=h,
\qbox{where $h\in K$ and
$w(x,y,\dots)$ is an element of the free group
$F(x,y,\dots)$,}
$$
solvable in $G$, is solvable in $K$.
It is known that, by a change of variables,
any such equation
can be transformed into the form
$$
x^mu(x,y,\dots)=h,
\hbox{ where $h\in K$ and
$u(x,y,\dots)$ lies in
the
commutator subgroup
of the
free group
$F(x,y,\dots)$}.
\eqno{(1)}
$$
Suppose that equation (1) has a solution
$(\~x,\~y,\dots)$
in $G$.
Since there is no principal
difference between $d_i$,
we can assume that
$$
\~x=d_1^\epsilon b^la_1^{k_1}a_2^{k_2}a_3^{k_3},
\qbox{where $\epsilon,l,k_i\in\Z$}.
\eqno{(2)}
$$
We have a homomorphism \emph{the first coordinate}:
$$
f\:G\to D_\infty=\gp{b'}_2\semitimes\gp{a'}_\infty,
\qbox{where
$f(a_1)=a'$,\quad $f(a_2)=f(a_3)=f(d_1)=1$,\quad $f(b)=f(d_2)=f(d_3)=b'$,
}
$$
and the natural homomorphism \emph{degree}\/:
$$
\deg\:G\to\Z,
\qbox{where
$\deg(a_i)=\deg(d_i)=0$,\quad $\deg(b)=1$,
}
$$
Applying these homomorphisms to the given equality
$\~x^mu(\~x,\~y,\dots)=h$,
we obtain
$$
f\Bigl(\~x^mu(\~x,\~y,\dots)\Bigr)=f(\~x)^mu(f(\~x),f(\~y),\dots)=f(h),
\qqbox{and}
\deg\Bigl(\~x^mu(\~x,\~y,\dots)\Bigr)=m\cdot\deg(\~x)=\deg(h).
\eqno{(3)}
$$
Now, consider the elements $\^x,\^y,\dots\in K$ obtained
from $\~x,\~y,\dots\in G$ by the
changes
$$
\eqalign{
a_1\mapsto a=a_1a_2a_3,
\quad
a_2\mapsto 1,
\quad
a_3\mapsto 1,
\quad
d_1\mapsto 1,
\quad
d_2\mapsto b,
\quad
d_3\mapsto b
\qbox{(and $b\mapsto b$)},&
\cr
\qbox{\sl which preserve the first coordinate of any element.
}&
}
\eqno{(4)}
$$
For instance, the element $\~x$, given by expression (2), turns into
$$
\^x=b^la^{k_1}=b^la_1^{k_1}a_2^{k_1}a_3^{k_1}.
\eqno{(5)}
$$
We claim that the tuple $(\^x,\^y,\dots)$ is a solution to equation (1) in
$K$. Indeed,
\-
the first coordinates of
$\^x^mu(\^x,\^y,\dots)$ and $h$
are the same:
\newline
$
f\Bigl(\^x^mu(\^x,\^y,\dots)\Bigr)
=
f(\^x)^mu(f(\^x),f(\^y),\dots)
\=^4
f(\~x)^mu(f(\~x),f(\~y),\dots)
\=^3
f(h)
$
(where (4) and (3) imply the corresponding equalities);
\-
and the degrees of
$\^x^mu(\^x,\^y,\dots)$ and $h$
are the same:
$
\deg\Bigl(\^x^mu(\^x,\^y,\dots)\Bigr)=
m\cdot\deg(\^x)
\=^5
ml
\=^2
m\cdot\deg(\~x)
\=^3
\deg(h),
$

\enditem
It remains to note that an element of $K\subset G$
is uniquely determined by its first coordinate and degree.
Thus, $\^x^mu(\^x,\^y,\dots)=h$, we have found a solution to equation
(1) in $K$, and this completes the proof.

\s 2.
Two lemmata on quotient groups

\proclaim
Lawfication lemma
{\rm([RKh13], Lemma 1.1)}.
If $V(G)$ is a verbal subgroup of a group $G$, and $H$ is
a verbally
closed subgroup of $G$, then
$H\cap V(G)=V(H)$ \(i.e. the verbal subgroup of $H$
corresponding to the same variety\),
and
the image $H/V(H)\subseteq G/V(G)$
of $H$ under the natural
homomorphism $G\to G/V(G)$ is verbally closed in $G/V(G)$.

\proclaim
Dihedral-quotient lemma.
If the Klein bottle
group
$K=\pres<a,b|a^b=a^{-1}>$ is a verbally closed subgroup
of a
group~$G$, then
\item{\rm 1)}
$\nc{b^2}\cap K=\gp{b^2}$, where
$\nc{b^2}$ is the normal closure of $b^2$
in $G$;
\item{\rm 2)}
the subgroup $D_\infty=K/\gp{b^2}\subseteq G/\nc{b^2}$
is
verbally closed in $G/\nc{b^2}$.

\Proof
We can assume that $G$
satisfies the law $[x^2,y^2]=1$,
(because this law holds in $K$, and, therefore,
the lawfication lemma allows us to replace
$G$ with its quotient group by the corresponding verbal subgroup).

For groups with such a law,
as well as for all metabelian groups,
assertion 1) is a general
fact:
\disp{\sl
if $A$ is a normal abelian subgroup of a group $G$,
and the quotient group $G/A$ is also abelian,
then the intersection $C_A(X)$ of $A$ and
the centraliser of any subset $X\subseteq G$ is
normal in $G$.
}%
Indeed,
$\bigl(C_A(X)\bigr)^g=C_A(X^g)\supseteq C_A(XA)=C_A(X)$.

To obtain assertion 1), we put
$A=\gp{\{g^2\;|\;g\in G\}}$ and $X=K$;
we even get more than 1):
$$
\hbox{\sl the subgroups
$\nc{b^2}$ and $K$ commute in $G$
}.
$$

Let us prove 2) now.
Suppose that an
equation
$$
w(x,y,\dots)=h\nc{b^2},
\qbox{where $h\in K$ and $w(x,y,\dots)\in F(x,y,\dots)$,}
\eqno{(*)}
$$
is solvable in $G/\nc{b^2}$.
We have to show that this equation is solvable
in $D_\infty=K/\gp{b^2}\subseteq G/\nc{b^2}$.

\Case 1:
$h=1$.
In this case, equation $(*)$
has the trivial solution $(1,1,\dots)$
in $D_\infty$.

\Case 2:
$h=b$.
In this case,
the exponent sum of one of unknowns (say, $x$)
in the word $w$ has to be odd,
because otherwise the solvability
of equation $(*)$ in $G/\nc{b^2}$
would imply
a decomposition of $b$ into a product of squares in $G$,
that contradicts the verbal closedness of $K$
(because $b$ is not a product of squares in $K$,
even modulo $\gp a$).
An equation with odd exponent sum of $x$
and the right-hand side $b$
has in $D_\infty$ an obvious solution: $x=b,\ y=1,\ z=1,\dots$.

\Case 3:
$h=ba^k$.
This case is
the same as the preceding one, because
the
dihedral group has an automorphism mapping $ba^k$ to $b$.

\Case 4
{\rm(the last case modulo $\gp{b^2}$)}:
$h=a^k$, where $k\ne0$.
Suppose that
$\bigl(\~x\nc{b^2},\;\~y\nc{b^2},\dots\bigr)$ is
a solution to equation~$(*)$
in $G/\nc{b^2}$.
Then the equation
$$
[t,\bigl(w(x,y,\dots)\bigr)^2]=[b,h^2]=a^{4k}
\qbox{(where $t$ is a new unknown)}
$$
has a solution $(\~t=b,\;\~x,\;\~y,\dots)$  in $G$,
because $\nc{b^2}$ commutes with
$K$ as noted above.
The verbal closedness
of $K$ in~$G$
implies that this equation has
a solution in $K$, i.e.
$$
[\^t,\bigl(w(\^x,\^y,\dots)\bigr)^2]=a^{4k}
\qbox{for some $\^t,\^x,\^y,\dots\in K$}.
\eqno{(**)}
$$
This means that

\-
$\^t\in b\gp{a,b^2}$,
because all other elements of $K$
commute with squares;
we can even assume that $\^t=b$, since
$a$~and~$b^2$ commute with all squares and
do
not affect commutator~$(**)$;

\-
$\bigl(w(\^x,\^y,\dots)\bigr)^2\in \gp{b^2}\cup \gp{a^2,b^4}$,
because only these elements are squares in $K$;

\-
$\bigl(w(\^x,\^y,\dots)\bigr)^2\in a^{2k}\gp{b^4}$,
because, only for such elements
of $\gp{b^2}\cup \gp{a^2,b^4}$, the commutator
with $\^t=b$ gives~$a^{4k}$;

\-
$w(\^x,\^y,\dots)\in a^k\gp{b^2}$,
because each element of the coset
$a^{2k}\gp{b^4}$ has a unique square root
in $K$.

\enditem
We have found a solution $(\^x\gp{b^2},\^y\gp{b^2},\dots)$
to equation $(*)$ in $D_\infty=K/\gp{b^2}$;
this completes the proof.

\s 3.
Proof of the second assertion of the theorem

Suppose that
the fundamental group $K$ of the
Klein bottle is
a verbally closed subgroup
of a
finitely generated group~$G$.
We have to construct a retraction onto $K$
from an index-at-most-two
subgroup of $G$
containing $K$.

The group $G$ has two normal subgroups:
$$
N_1=G'
\hbox{ --- the commutator subgroup, and }
N_2=\nc{b^2}
\hbox{ --- the normal closure of $b^2$}.
$$
Taking the quotients transforms $K$
into
$$
K/K'=\gp a_2\times\gp b_\infty\subseteq G/N_1
\hbox{ and }
K/\gp{b^2}=D_\infty\subseteq G/N_2
$$
(by the lawfication
lemma
and
dihedral-quotient lemma); these images
of $K$ in $G/N_i$ remain verbally closed in $G/N_i$
(by the same lemmata).
Therefore, $K/K'$ and~$K/\gp{b^2}$ are retracts
of $G/N_i$,
because all abelian groups [Mazh18]
and the infinite dihedral group [KMM18]
are
strongly verbally closed.

Thus,
we obtain the epimorphisms
$$
\deg\:G\to G/G'\to K/K'=\gp a_2\times\gp b_\infty\to\gp b_\infty
\too^\iso\Z
\qqbox{and}
f\:G\to G/\nc{b^2}\to K/\gp{b^2}=D_\infty
$$
such that $\deg(b)=1$, $f(b)=b\gp{b^2}$, $f(a)=a\gp{b^2}$.

Combining these two ``pseudo-retractions", we construct
the homomorphism
$$
\Phi\:G\to \Z\times D_\infty,
\quad
g\mapsto\bigl(\deg(g),\;f(g)\bigr).
$$
The restriction $\phi$ of $\Phi$ to
$K$ is injective, and the image of this restriction is
so-called \emph{fibered product}\/:
$$
\phi(K)=
\Phi(K)=\Biggl\{\Bigl(i,\;b^ja^k\gp{b^2}\Bigr)\;\Biggm|\;
i\equiv j\pmod2\Biggr\}
\hbox{
--- an index-two subgroup of $\Z\times D_\infty$.
}
$$
Therefore, the subgroup $\Phi^{-1}(\Phi(K))\subseteq G$ has index
at most two in $G$ and admits a retraction onto $K$:
$$
G\supseteq\Phi^{-1}(\Phi(K))\too^\Phi\Phi(K)\too^{\phi^{-1}}K.
$$

\References

[KlMa18]
A. A. Klyachko, A. M. Mazhuga,
Verbally closed virtually free subgroups,
Sbornik: Mathematics, 209:6 (2018), 850-856.
\arXiv 1702.07761

[Mazh19]
A. M. Mazhuga,
Free products of groups are strongly verbally closed,
Sbornik: Mathematics, 210:10 (2019), 1456-1492.
\arXiv 1803.10634

[RT19]
V. A. Roman'kov, E. I. Timoshenko,
On verbally closed subgroups of free solvable groups,
Vestnik Omskogo Universiteta, 24:1 (2019), 9-16
(in Russian).
\arXiv 1906.11689

[RKh13]
V. A. Roman'kov  and N. G. Khisamiev.
Verbally and existentially closed subgroups of free nilpotent groups,
Algebra and Logic, {52:4} (2013), 336-351.

[BMR99]
G. Baumslag, A. Myasnikov, V. Remeslennikov,
Algebraic geometry over groups I. Algebraic sets and ideal theory,
J. Algebra, 219:1 (1999), 16-79.

[Bog18]
O. Bogopolski,
Equations in acylindrically hyperbolic groups and verbal closedness,
arXiv:1805.08071.

[Bog19]
O. Bogopolski,
On finite systems of equations in acylindrically hyperbolic groups,
arXiv:1903.10906.

[KMM18]
A. A. Klyachko, A. M. Mazhuga, V. Yu. Miroshnichenko,
Virtually free finite-normal-subgroup-free groups
are strongly verbally closed,
J. Algebra, 510 (2018), 319-330.
\arXiv 1712.03406

[Lee02]
D. Lee,
On certain C-test words for free groups,
J. Algebra, 247:2 (2002), 509-540.
\arXiv math/0103108

[Mazh17]
A. M. Mazhuga,
On free decompositions of verbally closed subgroups
of free products of finite groups,
J.~Group Theory, 20:5 (2017), 971-986.
\arXiv 1605.01766

[Mazh18]
A. M. Mazhuga,
Strongly verbally closed groups,
J. Algebra, 493 (2018), 171-184.
\arXiv 1707.02464

[MR14]
A. Myasnikov, V. Roman'kov,
Verbally closed subgroups of free groups,
J. Group Theory, 17:1 (2014), 29-40.
\arXiv 1201.0497

[New85]
M. Newman,
A note on Fuchsian groups,
Illinois J. Math., 29:4 (1985), 682-686.

[Rom12]
V. A. Roman'kov,
Equations over groups,
Groups Complexity Cryptology,
4:2 (2012), 191-239.

\end